\documentclass[a4paper,10pt]{amsart}
\usepackage{amsmath, amsthm, amssymb, mathrsfs}
\usepackage{graphicx}
\usepackage{verbatim}
\usepackage[hang, small, bf]{caption}
\usepackage{hanging}
\usepackage{multirow}
\usepackage{enumerate}
\usepackage[all, cmtip]{xy}
\usepackage{fullpage}
\usepackage{cases}

\theoremstyle{remark}
\newtheorem{thm}{Theorem}[section]
\newtheorem{lemma}[thm]{Lemma}

\newtheorem{conj}[thm]{Conjecture}
\newtheorem{rem*}[thm]{Remark}
\newtheorem*{thmA}{Theorem A}
\newtheoremstyle{example}{\topsep}{\topsep}%
     {}%         Body font
     {}%         Indent amount (empty = no indent, \parindent = para indent)
     {\bfseries}% Thm head font
     {}%        Punctuation after thm head
     {\newline}%     Space after thm head (\newline = linebreak)
     {\thmname{#1}\thmnumber{ #2}\thmnote{ #3}}%         Thm head spec
\newtheorem*{sethm}{Theorem~\ref{sethm}}

\theoremstyle{example}
\newtheorem{example}[thm]{Example}

\theoremstyle{definition}
\newtheorem{Def}{Definition}[section]

\newcommand{\ZZ}{\mathbb{Z}}

\newcommand{\Cspace}[2]{\mathcal{C}^{#1}({#2}) }
\newcommand{\UCspace}[2]{U\mathcal{C}^{#1}({#2}) }
\newcommand{\Dspace}[2]{\mathcal{D}^{#1}({#2}) }
\newcommand{\UDspace}[2]{U\mathcal{D}^{#1}({#2}) }
\newcommand{\VT}{VT} % Vertex terminal
\newcommand{\ET}{ET} % Edge terminal
\newcommand{\goesto}{\rightarrow}
\newcommand{\into}{\hookrightarrow}
\newcommand{\collapseto}{\searrow}

\DeclareMathOperator{\rank}{rank}

\begin{document}
\title{Abrams's Stable Equivalence for Graph Braid Groups}

\author[ucd]{Paul Prue}
%\ead{pruep@scc.losrios.edu}

\author[ucd]{Travis Scrimshaw}
%\ead{tscrim@ucdavis.edu}

%\address{University of California, Davis, One Shields Avenue, Davis, CA 95616}

%\maketitle

\begin{abstract}
In his PhD thesis~\cite{abramsthesis}, Abrams proved that, for a natural number $n$ and a graph $G$ with at least $n$ vertices, the $n$-strand configuration space of $G$, denoted $\Cspace{n}{G}$, deformation retracts to a compact subspace, the \emph{discretized} $n$-strand configuration space, provided $G$ satisfies two conditions: each path between distinct \emph{essential vertices} (vertices of degree not equal to 2) is of length at least $n+1$ edges, and each path from a vertex to itself which is not nullhomotopic is of length at least $n+1$ edges. Using Forman's discrete Morse theory for CW-complexes, we show the first condition can be relaxed to require only that each path between distinct essential vertices is of length at least $n-1$.
\end{abstract}


\maketitle

%%%%%%%%%%%%%%%%%%%%%%%%%%%%%%%%%%%%%%%
\section{Introduction}
\label{sec_intro}
%%%%%%%%%%%%%%%%%%%%%%%%%%%%%%%%%%%%%%%

The goal of this paper is to establish sufficient conditions such that a \emph{braid group on a graph} may be studied via a certain CW complex associated to the graph. Let $X$ denote a connected topological space. An \emph{$n$-strand configuration in $X$} is an $n$-point subset of $X$. The \emph{unordered $n$-strand configuration space} of $X$ is the space of unordered subsets consisting of $n$ distinct elements of $X$, and is denoted $\UCspace{n}{X}$. (We use the terms \emph{labeled} and \emph{unlabeled} as synonyms for the terms \emph{ordered} and \emph{unordered}, respectively.) For a positive integer $n$, the classical braid group $B_n$ is the fundamental group $\pi_1(\UCspace{n}{D^2})$, where $D^2$ is the 2-dimensional topological disk. Thus, from the configuration-space perspective, a braid is simply a loop in the space $\UCspace{n}{D^2}$. Similarly, the \emph{ordered $n$-strand configuration space} of $D^2$, denoted $\Cspace{n}{D^2}$, is the space of ordered tuples consisting of $n$ distinct elements of $X$. The classical \emph{$n$-strand pure braid group}, denoted $PB_n$, is the fundamental group of the ordered $n$-strand configuration space of a disk. Note, the quotient map from the configuration space $\Cspace{n}{D^2}$ to the unordered configuration space $\UCspace{n}{D^2}$ induces a short exact sequence $1 \goesto PB_n \goesto B_n \goesto \Sigma_n \goesto 1$, where $\Sigma_n$ is the symmetric group on $n$ symbols. For an extensive reference on classical braid groups, see~\cite{birman-2004}. 

In the case of graph braid groups, we let $X = G$ be a connected graph, viewed as a 1-dimensional CW complex. The ordered $n$-strand configuration space of $G$ is denoted $\Cspace{n}{G}$. The \emph{$n$-strand pure graph braid group} $PB_n(G)$ is the fundamental group $\pi_1(\Cspace{n}{G})$. The unordered configuration space is $\UCspace{n}{G}$, and its fundamental group $\pi_1(\UCspace{n}{G})$ is the \emph{$n$-strand graph braid group} $B_n(G)$. Note that these fundamental groups do not depend on basepoint. Usually the configuration space is connected, but even when it is disconnected the components are all homeomorphic~\cite{abramsthesis}.

\begin{figure}[!t]
 \begin{center}
 \includegraphics[scale=0.5]{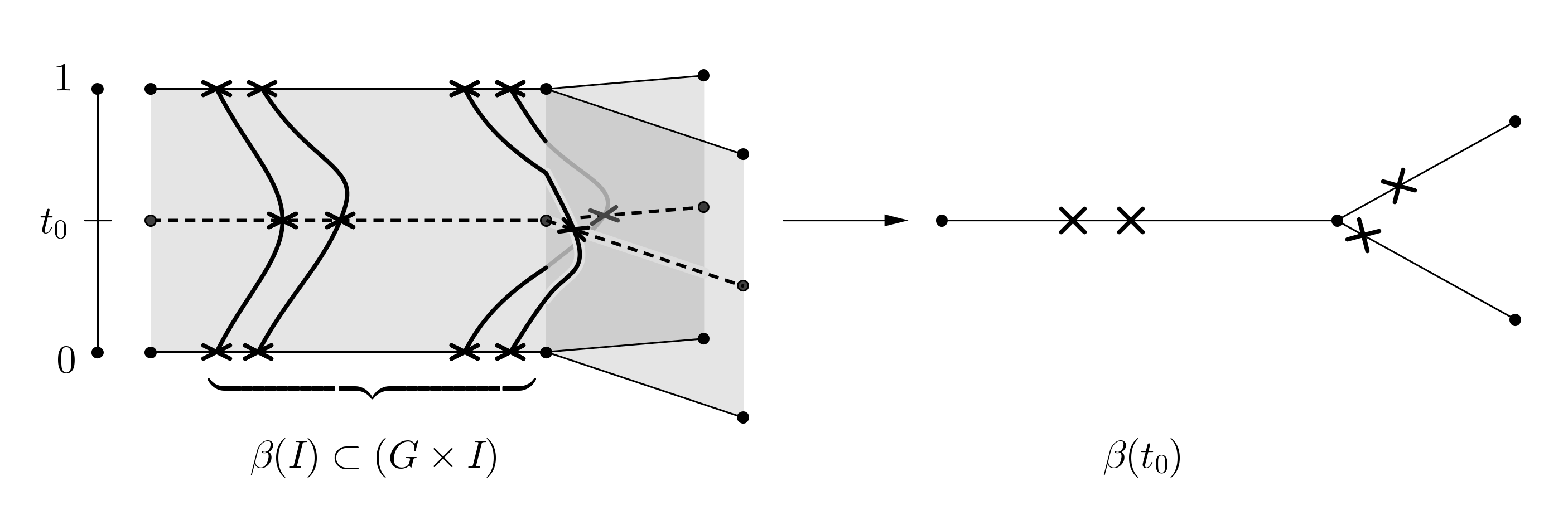}
 \end{center}
 \caption{A nontrivial 4-braid in the cylinder $K_{3,1} \times [0,1]$. At each $t \in [0,1]$, an $n$-braid defines a configuration of $n$ points of the graph, illustrated by the \sf X\rm's on the graph at right.}
 \label{graphcyl}
\end{figure}

Graph braid groups, like classical braid groups, can also be viewed as isotopy classes rel endpoints of \emph{braids} (i.e., certain $n$-tuples of pairwise-disjoint paths) in the cylinder on a topological space. For classical braid groups, this cylinder is $D^2 \times I$, where $I$ is the interval $[0,1]$. In the case of the graph braid group $B_n(G)$, one considers instead braids in the cylinder $G \times I$. Figure~\ref{graphcyl} shows a non-trivial braid in $G \times I$, where $G$ is isomorphic to the complete bipartite graph $K_{3,1}$. A braid $\beta : I \goesto G \times I$ can be thought of as describing the simultaneous and continuous movements of the $n$ strands, or \emph{tokens}, without collisions, on $G$. To each $t \in I$, the map $\beta$ associates a configuration $\beta(t)$ of the $n$ strands on $G$. Since $\beta$ is a loop in the (ordered/unordered) configuration space, it follows that the configurations $\beta(0)$ and $\beta(1)$ are equal. For example, in Figure~\ref{graphcyl}, $\beta(0) = \beta(1)$ as unordered configurations.

Besides providing a class of interesting mathematical objects, graph braid groups have real-world applications that have been discussed in~\cite{abramsghristfinding} and \cite{ghrist-1999}. An example often given is that of a fleet of mobile robots inside a factory, whose movement is confined to a network of track or guide tape. For an idealized robot of infinitesimal size, the configuration space of points on a graph shaped like the track network exactly describes the space of configurations of the fleet of robots. 

Abrams introduces the notion of a \emph{discretized} configuration space of $n$ strands on a graph $G$ in his PhD thesis~\cite{abramsthesis}. This is a compact subspace of a configuration space of $G$, consisting of only those $n$-strand configurations $x$ in which, for each pair of strands in $x$, every path in $G$ between the two strands contains at least one full edge of $G$. Note that $\Cspace{n}{G}$ is a subspace of the cubical complex $\prod^n G$, but does not inherit its CW structure, as it is not a compact subspace. In contrast, the \emph{discretized labeled configuration space} of a graph $G$, denoted $\Dspace{n}{G}$, has a CW complex structure as a subcomplex of $\prod^n G$, as does the \emph{discretized unlabeled configuration space}, denoted $\UDspace{n}{G}$. The space $\UDspace{n}{G}$ is obtained as a quotient of $\Dspace{n}{G}$ by the action of the symmetric group $\Sigma_n$, which acts by permuting the coordinates of a labeled configuration. Interesting examples among these discretized configuration spaces have been described by Abrams and Ghrist (see~\cite{abramsthesis}, \cite{abramsghristfinding}, and \cite{ghrist-1999}). We include an example (\ref{collapse-example}) in this paper illustrating the discrete Morse function defined in the proof of Lemma~\ref{sclemma}.

For a given $n$, Theorem 2.1 of~\cite{abramsthesis} gives sufficient conditions on $G$ to guarantee that $\Cspace{n}{G}$ deformation retracts onto the subspace $\Dspace{n}{G}$, a cubical complex. We state the theorem here for reference. An \emph{essential vertex} of a graph is a vertex whose degree is not equal to 2.
\begin{thmA} \cite{abramsthesis}
Let $G$ be a graph with at least $n$ vertices. Then $\Cspace{n}{G}$ deformation retracts onto $\Dspace{n}{G}$ if
\begin{enumerate}[(A$^\prime$)]
\item each path connecting distinct essential vertices of $G$ has length at least $n+1$, and
\item each homotopically essential path connecting a vertex to itself has length at least $n+1$.
\end{enumerate}
\label{abramsse}
\end{thmA}
The homotopy equivalence established by Theorem~A guarantees that one obtains the $n$-strand pure braid group of $G$ as $\pi_1(\Dspace{n}{G})$, provided $G$ is subdivided as in the hypotheses of the theorem. However, condition (A$^\prime$) of Abrams's theorem is not best possible, as he conjectures.  Several papers (including~\cite{abramsghristfinding}, \cite{farley-sabalka-2005}, and \cite{sabalkathesis}) have cited the theorem incorrectly, assuming the improved bounds on subdivision for which we give proof in Theorem~\ref{sethm}.

In this paper, we prove that the homotopy type of $\Dspace{n}{G}$ stabilizes, given sufficient subdivision on the edges of $G$. We use Forman's discrete Morse theory to prove this statement, Lemma~\ref{sclemma}, giving an explicit discrete Morse function on $\Dspace{n}{G}$ under the assumption that $G$ has ``excess'' subdivision on (at least) one edge. Having established Lemma~\ref{sclemma}, we use a straightforward direct limit argument to prove the sufficiency of the improved bounds on graph subdivision, conjectured by Abrams, such that the $n$-strand configuration space of a graph $G$ deformation retracts to the $n$-strand discretized configuration space. This is the Stable Equivalence Theorem stated below, and proved in Section~\ref{main}.

\begin{sethm}[Stable Equivalence]
Let $n>1$ be an integer, and $G$ a finite, connected graph with at least $n$ vertices. The $n$-point configuration space $\Cspace{n}{G}$ deformation retracts onto the discretized configuration space $\Dspace{n}{G}$ if
\begin{enumerate}[(A)]
%---------------(A)--------------------
\item each path connecting distinct essential vertices of $G$ has length at least $n-1$, and
%---------------(B)--------------------
\item each homotopically essential path connecting a vertex to itself has length at least $n+1$.
\end{enumerate}
\end{sethm}

The result in the theorem extends also to the unlabeled configuration spaces, giving sufficient conditions on subdivision in $G$, for a given $n$, such that $\UCspace{n}{G}$ deformation retracts onto $\UDspace{n}{G}$. We note that when $n=1$, both the topological and discretized configuration spaces are homeomorphic to $G$.

We repeat the conjecture, made by Abrams in~\cite{abramsthesis}, that conditions (A) and (B) of the theorem are in fact necessary (hence optimal).

\begin{conj}
If $n > 1$ is an integer and $G$ is a finite, connected graph with at least $n$ vertices, the space $\Cspace{n}{G}$ deformation retracts onto $\Dspace{n}{G}$ if \emph{and only if} $G$ and $n$ satisfy conditions (A) and (B) of Theorem~\ref{sethm}.
\end{conj}
On the basis of Theorem~\ref{sethm} and the conjecture above, we make the following provisional definition.

\begin{Def}
Let $G$ be a finite, connected graph, and let the integer $n > 1$ be given. We say $G$ is \emph{sufficiently subdivided for $n$ strands} if $G$ satisfies conditions (A) and (B) of Theorem~\ref{sethm}.
\label{sufficientsubdivisiondef}
\end{Def}
Where $n$ is clear from the context, we simply say $G$ is sufficiently subdivided.

Many results about graph braid groups are known. In~\cite{ghrist-1999}, Ghrist showed that the unordered configuration spaces $\UCspace{n}{G}$ are Eilenberg-Maclane spaces of type $K(B_n(G),1)$. Abrams showed in~\cite{abramsthesis} that graph braid groups have solvable word and conjugacy problems, as they are fundamental groups of locally CAT(0) cubical complexes. Farley and Sabalka use Forman's discrete Morse theory in~\cite{farley-sabalka-2005} and~\cite{farley-sabalka-2012} to calculate presentations for graph braid groups. Farley showed in~\cite{farley-2006} that braid groups of trees have free abelian integral homology groups in every dimension, and in~\cite{farley-2007} computed presentations of the integral cohomology rings $H^{*}(\UCspace{n}{G}; \ZZ)$, where $T$ is any tree and $n$ is arbitrary. Farley and Sabalka proved in~\cite{farley-sabalka-2008} that a graph braid group of a tree is right-angled Artin if the tree is linear or $n < 4$. Crisp and Wiest showed in~\cite{crisp-wiest} that every graph braid group embeds in a right-angled Artin group, hence, graph braid groups are linear, bi-orderable, and residually finite.

Abrams describes an explicit deformation retraction $\Cspace{n}{G} \goesto \Dspace{n}{G}$ when $n=2$ in~\cite{abramsthesis}. In~\cite{kimkopark}, Kim, Ko, and Park prove a version of Theorem~\ref{sethm} under the assumption that $n\ge 3$. Their proof depends on Abrams's thesis, and employs the discrete Morse-theoretic tools of~\cite{farley-sabalka-2005}. In contrast, we prove the theorem using an explicitly defined discrete Morse function, giving a unified treatment for all indices $n>1$.

The outline of this paper is as follows. In Section~\ref{lemmas}, we use Forman's discrete Morse theory -- in particular, we explicitly construct a discrete Morse function -- to prove that the homotopy type of a discretized configuration space eventually stabilizes after repeated subdivision of a single edge. Section~\ref{main} contains our main result and an example illustrating the discrete Morse technique of Section~\ref{lemmas}.

The authors gratefully acknowledge Lucas Sabalka for his mentorship and encouragement. This paper grew out of a summer 2008 REU with him at UC Davis, supported by NSF VIGRE grant DMS-0636297.

%%%%%%%%%%%%%%%%%%%%%
\section{Discrete Morse Function}
\label{lemmas}
%%%%%%%%%%%%%%%%%%%%%

The proof of Lemma~\ref{sclemma} uses Forman's discrete Morse theory, so we review some of the basic notions of that theory here.

Let $X$ be a finite regular CW complex. The notation $\alpha^{(r)}\in X$ means that $\alpha$ is a cell of $X$ whose dimension is $r$. If
$\gamma^{(r-1)} \subset \overline{\alpha^{(r)}}$ and $\gamma\ne\alpha$, we say $\gamma$ is a \emph{face} of $\alpha$, and write $\gamma < \alpha$. A function $f \colon X \to \mathbb{R}$ is a \emph{discrete Morse function} if, for every $\alpha^{\left(r\right)}\in X$, both
\[
L_{\alpha}:=\#\left\{ \gamma^{\left(r-1\right)}<\alpha\mid f\left(\gamma\right)\ge f\left(\alpha\right)\right\} \le 1
\]
and
\[
U_{\alpha}:=\#\left\{ \beta^{\left(r+1\right)}>\alpha\mid f\left(\beta\right)\le f\left(\alpha\right)\right\} \le 1.
\]

Given a discrete Morse function on $X$, we classify cells of $X$ as follows. If $L_{\alpha}=0$ and $U_{\alpha} = 1$, then $\alpha$ is \emph{redundant}. If $L_{\alpha}=1$ and $U_{\alpha}=0$, then $\alpha$ is \emph{collapsible}. If $L_{\alpha}=U_{\alpha}=0$, then $\alpha$ is \emph{critical}. (By Lemma 2.5 of~\cite{forman-1998}, $U_{\alpha}$ and $L_{\alpha}$ cannot both be positive when $f$ is Morse.)

In the language of discrete gradient vector fields~\cite{forman-2002}, indicated by ``arrows,'' a redundant cell is the tail of an arrow (it is in the domain of the gradient vector field), while a collapsible cell is the head of an arrow (in the range of the vector field). A critical cell is neither the head nor the tail of an arrow.

Given a cell complex $X$ with a discrete Morse function $f$, for any real number $p$, the \emph{level subcomplex} $M(p)$ is defined by
\[
 M(p) = \bigcup_{\substack{\alpha\in X \\ f(\alpha)\le p} } \left( \bigcup_{ \gamma < \alpha} \gamma \right).
\] 
(In general, ``$\gamma < \alpha$'' indicates that $\gamma$ is a face of $\alpha$ of any codimension. For our purposes, however, it suffices to assume codimension 1.) 

We will also need the following theorem of Forman. To indicate the cellular collapse from a CW-complex $X$ onto a subcomplex $Y$, we write $X \collapseto Y$.
\begin{thm}[Thm.~3.3,~\cite{forman-1998}]
\label{thm:collapse_level}
If $a < b$ are real numbers such that the interval $[a, b]$ contains no critical values of a discrete Morse function $f$, then $M(b) \collapseto M(a)$.
\end{thm}
Applying Theorem~\ref{thm:collapse_level} to the discrete Morse function $f$ defined in Equation~\eqref{eq:morse} below, we prove the following lemma, which shows a simple homotopy equivalence between certain level subcomplexes of a given discretized graph configuration space.

\begin{lemma}
\label{sclemma}
Let $G$ be a (finite, connected) graph and $n$ an integer greater than 1. Fix an integer $m \geq n-1$. Suppose $p$ is a path in $G$ consisting of the $m$ edges $\left\{ v_{0}, v_{1} \right\}, \left\{ v_{1},v_{2}\right\}, \dotsc, \left\{ v_{m-1}, v_{m}\right\}$, where
\begin{itemize}
\item $v_{i} \neq v_{j}$ if $i \neq j$, unless $m \geq n+1$ and $\left\{ i,j \right\} = \left\{ 0,m \right\}$;
and
\item $\deg\left(v_{i}\right) = 2$ for $1 \leq i \leq m-1$.
\end{itemize}

Let $G^{\prime}$ be the graph obtained from $G$ by subdividing some edge of $p$. Then $\Dspace{n}{G^{\prime}} \collapseto \Dspace{n}{G}$.
\end{lemma}

\begin{rem*}
In fact, $\Dspace{n}{G}$ is not a true subcomplex of $X = \Dspace{n}{G^{\prime}}$, but there is a subcomplex $Y\subset X$ that corresponds to the image of $\Dspace{n}{G}$ under inclusion (namely, the cells of $X$ of zero rank, as described in the proof).
\label{not-true-subcomplex}
\end{rem*}

\begin{proof}
We use the following notation. Let $x = \left(x_{1},x_{2},\ldots,x_{n}\right)$ be a cell of $\Dspace{n}{G}$, with each $x_{i}$ a cell (vertex or open edge) of $G$, and $\overline{x_i} \cap \overline{x_j} = \emptyset$ for $i\ne j$. By $\overline{x}$ we denote the subgraph of $G$ consisting of the edges and vertices in $x$:
\[
\overline{x} = \bigcup_{i=1}^{n} \overline{x_{i}}.
\]
Thus, for a vertex $v$ of $G$, the notation ``$v \in \overline{x}$'' means that $v$ is either a vertex or the endpoint of an edge in $x$. The notation ``$v \in x$'' means that the vertex $v$ is a factor of $x$ (i.e., the cell $x_{i} = v$ for some $i$), and consequently, \emph{no edge} with endpoint $v$ is in $x$. Similarly, for an edge $e$, ``$e \in x$'' means $x_i = e$ for some $i$.

First, we label the vertices of $p^{\prime} \subseteq G^{\prime}$,
from $v_{0}$ through $v_{m+1}$, and suppose $v_{i}$ is the added
vertex for some $1 \leq i \leq m$. We label the edges of $p^{\prime}$
as follows.
\begin{align*}
e_{1} = & \left\{ v_{i-1},v_{i}\right\} \\
e_{2} = & \left\{ v_{i-2},v_{i-1}\right\} \\
\vdots & \\
e_{i} = & \left\{ v_{0},v_{1}\right\} \\
e_{i+1} = & \left\{ v_{i},v_{i+1}\right\} \\
\vdots & \\
e_{m+1} = & \left\{ v_{m},v_{m+1}\right\}
\end{align*}
The endpoints of an edge $e \subset p^{\prime}$ are $\iota e$ and
$\tau e$, such that there is a geodesic segment $\overline{v_{i}\;\tau e}\subset p^{\prime}$
that contains $\iota e$ but not $v_{0}$. Note the repetition $\iota e_{1} = \iota e_{i+1} = v_{i}$.

The \emph{rank} of a cell $x$ of $\Dspace{n}{G^{\prime}}$
is defined as follows.
\begin{itemize}
\item
Set $\rank(x) = 0$ if less than two cells of $x$ meet the subgraph
$\overline{e_{1} \cup e_{i}}$, or if no cell of $x$ meets $v_{i}$.

\item Otherwise, $\rank(x)$ is the positive number
\[
\rank(x) = \min\left\{ j \mid x\cap\tau e_{j} = \emptyset \mbox{ or } e_{j}\in x \right\} .
\]
\end{itemize}

We call a cell of rank 0 \emph{remote}. The cells of positive rank
are classified as follows. We say $x$ is \emph{vertex terminal of
rank $q$}, and write $x\in \VT_{q}$, if $0< \rank(x) = q$ and $e_{q} \notin x$.
On the other hand, $x$ is \emph{edge terminal of rank $q$}, denoted
$x \in \ET_{q}$, if $0 < \rank(x) = q$ and $e_{q}\in x$.

Let $f$ be the discrete function defined on the cells of $\Dspace{n}{G^{\prime}}$
as follows. Suppose $x^{(r)}$ is a cell of $\Dspace{n}{G^{\prime}}$. We set
\begin{equation}
\label{eq:morse}
f \left(x\right) = \begin{cases}
r & \text{if } \rank(x) = 0, \\
r+q+n+1 & \text{if } x\in \VT_{q}, \\
r+q+n & \text{if } x\in \ET_{q}.
\end{cases}
\end{equation}
We will show in Lemma~\ref{f_morse-lemma} below that $f$ is a discrete Morse function, whose critical cells are precisely those of rank 0. Consequently, all critical values of $f$ lie in the interval $\left[0,n\right]$. By Theorem~\ref{thm:collapse_level}, the level subcomplex $M\left(3n+2\right) = \Dspace{n}{G^{\prime}}$ collapses to $M\left(n+1/2\right) = Y$, the subcomplex described in Remark~\ref{not-true-subcomplex}.
\end{proof}

\begin{lemma}
\label{f_morse-lemma}
The function $f$ defined in~\eqref{eq:morse} is a discrete Morse function.
\end{lemma}

\begin{proof}
The defining conditions for a Morse function are local, involving face/coface pairs. Thus, the following facts are relevant. Each is a consequence of the simple observation that, in the cubical complex $\Dspace{n}{G^\prime}$, a coface $y^{(r+1)}$ of a cell $x^{(r)}$ is obtained from $x$ by a substitution $v \in \partial e \mapsto e$ for some $v \in x$, $e \in y$. Conversely, a face $x^{(r)}$ of a cell $y^{(r+1)}$ comes from a substitution
$e \mapsto v \in \partial e$.
\begin{itemize}
\item Every face of a remote face is remote. Put another way, the set $R$
of remote faces is a subcomplex of $\Dspace{n}{G}$.
\item If $x \in \VT_{q}$, and $y>x$, then $\rank(y) \geq \rank(x)$. (None of
the edges $e_{k}$, for all $1 \leq k \leq q$, is in $x$, but all of their endpoints
meet $x$.) In particular, $y$ is not remote. On the other hand,
if $w < x$, then $\rank(w) \leq \rank(x)$. Moreover, $w$ is not edge terminal.
\item If $x \in \ET_{q}$, and $y > x$, then $\rank(y) = \rank(x)$, and $y$ is
edge terminal. If $w < x$, then $\rank(w) \leq \rank(x)$.
\end{itemize}

The diagram below summarizes some of this information, distinguishing
the possible respective classifications of a pair of cells $\alpha < \beta$.
We write $R$, $\VT$, and $\ET$, respectively, for the sets of remote,
vertex terminal, and edge terminal cells. If $A$ and $B$ are two
of these sets, there is an arrow $A \longrightarrow B$ if $\alpha < \beta$
for some pair of cells $\alpha \in A$ and $\beta \in B$.
\[
\xymatrix{R \ar[d] \ar[rd] \\
\VT \ar[r] & \ET
}
\]
We have suppressed self-loops at each node: with some exceptions,
a cell of a given type may have faces/cofaces of the same type. 

These constraints on face/coface pairs reduce the number of cases
to be checked in verifying that $f$ is Morse. We proceed now to this
analysis.

Let $x \in \Dspace{n}{G^{\prime}}$ be a cell of dimension $r$.
\begin{enumerate}
\item Suppose $x$ is remote, so $f(x) = r$.

\begin{enumerate}
\item If $w^{(r-1)}<x$, then $w$ is also remote, so $f\left(w\right) = r-1 < f\left(x\right)$.
Thus, $L_{x}=0$.
\item If $y^{(r+1)}>x$, then
\[
f(y) = \begin{cases}
r+1 & \mbox{if } y \mbox{ is remote,}\\
r+q+n+2 & \mbox{if } y \in \VT_{q} \mbox{ for some } 1\le q \le n,\\
r+q+n+1 & \mbox{if } y \in \ET_{q} \mbox{ for some } 1\le q \le n.
\end{cases}
\]
So $f\left(y\right) > f\left(x\right)$, meaning $U_{x} = 0$.
\end{enumerate}

Since $U_{x} = L_{x}=0$, we conclude $x$ is a critical cell of $f$.

\item Suppose $x\in \VT_{q}$ for some $1\leq q \leq n$. Then $f\left(x\right) = r+q+n+1$.

\begin{enumerate}
\item If $w^{(r-1)}<x$, then $w$ is not edge terminal, and $\rank(w) \leq \rank(x)$.
\[
f\left(w\right)= \begin{cases}
r-1 & \text{if $w$ is remote},\\
r+p+n & \text{if $w \in VT_{p}$ for some $p \leq q$.}
\end{cases}
\]
In either case, $f\left(w\right) < f\left(x\right)$, so $L_{x}=0$.
\item If $y^{(r+1)} > x$, then $y$ is not remote. Suppose $\rank(y) = s$, for
some $s \geq q$. Then
\[
f\left(y\right) = \begin{cases}
r+s+n+1 & \text{if } y \in \VT_{s},\\
r+s+n & \text{if } y \in \ET_{s}.
\end{cases}
\]
So $f\left(y\right) \geq f\left(x\right)$, with equality only if $y \in \ET_{q}$.
There is exactly one such coface of $x$, obtained from $x$ by the
substitution $\iota e_{q} \mapsto e_{q}$. Hence, $U_{x} = 1$.
\end{enumerate}

Since $L_{x} = 0$ and $U_{x} = 1$, $x$ is redundant.

\item Lastly, suppose $x \in \ET_{q}$ for some $1 \leq q \leq n$. Then $f\left(x\right) = r+q+n$.

\begin{enumerate}
\item If $w^{(r-1)} < x$, then $\rank(w) \leq \rank(x)$. Suppose $\rank(w) = p$,
for some $0 \leq p \leq q$. Then
\[
f\left(w\right) = \begin{cases}
r-1 & \text{if $w$ is remote,}\\
r+p+n & \mbox{if } w \in VT_{p},\\
r+p+n-1 & \mbox{if } w \in ET_{p}.
\end{cases}
\]
Thus, $f\left(w\right) \leq f\left(x\right)$, with equality only if
$w\in \VT_{q}$. There is exactly one such face of $x$, obtained from
$x$ by the substitution $e_{q} \mapsto \iota e_{q}$. So $L_{x}=1$.
\item If $y^{(r+1)} > x$, then $y \in \ET_{q}$ also, so $f\left(y\right) = r+q+n+2 > f\left(x\right)$.
Thus, $U_{x}=0$.
\end{enumerate}

Since $L_{x}=1$ and $U_{x}=0$, $x$ is collapsible.

\end{enumerate}

This analysis shows that $f$ is Morse. Furthermore, the critical
values of $f$ are in the interval $\left[0,n\right]$, while $\left[n+1/2,3n+2\right]$
contains no critical values of $f$. Applying Theorem 3.3 of~\cite{forman-1998}, we obtain the simple
homotopy equivalence $\Dspace{n}{G^{\prime}} = M\left(3n+2\right) \collapseto M\left(n+1/2\right) = Y$.
\end{proof}

%%%%%%%%%%%%%%%%%%%%%
\section{Main Result}
\label{main}
%%%%%%%%%%%%%%%%%%%%%

In this section, we prove our main result. Our proof is similar to the proof of Theorem~A given in~\cite{abramsthesis}. First, the following lemma establishes that $(\Cspace{n}{G}, \Dspace{n}{G})$ can be viewed as a CW-pair.

\begin{lemma}
\label{cw-structure}
Let $n>1$ be an integer, and suppose the graph $G$ is sufficiently subdivided for $n$ strands. Then there exists a CW structure on $\Cspace{n}{G}$ such that the inclusion $\iota : \Dspace{n}{G} \into \Cspace{n}{G}$ is a cellular map.
\end{lemma}

We omit the proof of this lemma, remarking only that one such cellular structure is partly induced on the complement $\Cspace{n}{G} \setminus \Dspace{n}{G}$ by a discrete collection of level sets of the \emph{diameter} function $\delta: \Cspace{n}{G} \goesto (0,\infty)\subset\mathbb{R}$, defined on configurations $x\in \Cspace{n}{G}$ by $\delta(x) = \min\{d(x_i, x_j) \mid i\ne j\}$, where $d$ denotes the distance in $G$ between two strands of a configuration.

\begin{thm}[Stable Equivalence]
\label{sethm}
Let $n > 1$ be an integer, and $G$ a finite, connected graph. The space $\Cspace{n}{G}$ deformation retracts onto $\Dspace{n}{G}$ if
\begin{enumerate}[(A)]
%---------------(A)--------------------
\item each path connecting distinct essential vertices of $G$ has length at least $n-1$, and
%---------------(B)--------------------
\item each homotopically essential path connecting a vertex to itself has length at least $n+1$.
\end{enumerate}
\end{thm}

\begin{proof}
Lemma~\ref{cw-structure} implies that $(\Cspace{n}{G}, \Dspace{n}{G})$ is a CW-pair. We show that the inclusion of the subcomplex induces isomorphisms on all homotopy groups.

Since $G$ is sufficiently subdivided, it contains at least $n$ vertices, so each connected component of $\Cspace{n}{G}$ contains a 0-cell of $\Dspace{n}{G}$. Thus, the inclusion $\iota \colon \Dspace{n}{G} \into \Cspace{n}{G}$ induces a surjection on $\pi_0$.

Let $c$ be a loop in $\Dspace{n}{G}$ which bounds a disk $D$ in $\Cspace{n}{G}$. Since $D$ is compact, there exists a subdivision $G^\prime$ of $G$ such that, for every $x\in D$, at least one full edge of $G^\prime$ separates each pair of strands in $x$. Since $\Dspace{n}{G^{\prime}}$ deformation retracts onto $\Dspace{n}{G}$ by Lemma~\ref{sclemma}, the boundary loop $c$ is nullhomotopic in $\Dspace{n}{G}$; thus, $\iota$ induces an injection on $\pi_1$. 

By a similar argument, $\iota$ induces a surjection on $\pi_1$, and an injection on $\pi_0$. For, if conditions (A) and (B) hold, then every path $\alpha: ([0,1], \partial [0,1]) \to (\Cspace{n}{G}, \Dspace{n}{G})$ can be homotoped rel boundary to lie completely within $\Dspace{n}{G}$.

Lastly, $\Cspace{n}{G}$ and $\Dspace{n}{G}$ are aspherical~\cite{abramsthesis, ghrist-1999}, so the map $\iota$ induces bijections on $\pi_n$ for all $n \geq 2$. Therefore, by the Whitehead theorem, $\Dspace{n}{G}$ is a deformation retract of $\Cspace{n}{G}$.
\end{proof}

The subdivision in a graph satisfying the conditions in Definition~\ref{sufficientsubdivisiondef} is \emph{sufficient} in the following sense. For any connected graph $G$, adding a finite number of non-essential vertices to $G$ (that is, subdividing an edge $e$ of $G$ by removing the interior of $e$, and adding a vertex of degree 2 connected by an edge to each endpoint of $e$) gives a new graph $G^\prime$ which is homeomorphic to $G$. Thus, for fixed $n$, the space $\Cspace{n}{G}$ is homeomorphic to $\Cspace{n}{G^\prime}$. The deformation retraction whose existence is proved in Theorem~\ref{sethm} projects to the quotient $\UCspace{n}{G} \goesto \UDspace{n}{G}$, so the theorem holds for unordered spaces also. Although $\UDspace{n}{G}$ and $\UDspace{n}{G^\prime}$ may be very different spaces in a combinatorial sense, the theorem guarantees that if $\pi_1(\UDspace{n}{G})$ is isomorphic to $B_n(G)$, then $\pi_1(\UDspace{n}{G^\prime})$ is isomorphic to $B_n(G)$, too.

The following example illustrates the homotopy equivalence between $\UDspace{2}{G}$ and $\UDspace{2}{G^\prime}$ for a particular sufficiently subdivided graph $G$ and the graph $G^\prime$ obtained from $G$ by subdividing one edge.

\vspace{12pt}
\begin{example}
Let $G$ be the graph homeomorphic to the letter \textsf{P} illustrated at left in Figure~\ref{p-optimal}, with vertex set $V=\{1, 2, 3, 4\}$ and edge set $E=\left\{ \,\{1,2\},\{2,3\},\{3,4\},\{2,4\}\,\right\}$. The unlabeled discretized space of 2-strand configurations on $G$ is the space $\UDspace{2}{G}$ illustrated at right in Figure~\ref{p-optimal}. The number of cells of each dimension is as follows.
\begin{itemize}
 \item 0-cells: $\binom{|V|}{2} = \binom{4}{2} = 6$, each an unordered pair of vertices of $G$.
 \item 1-cells: $8$, each of the form $[e, v]$, where $e\in E$ is a closed edge, and $v\in V$ is a vertex, with $v \notin \overline{e}$.
 \item 2-cells: one, namely, $[\{1,2\}, \{3,4\}]$ -- since these are the only edges $e,e^\prime\in E$ with $\overline{e} \cap \overline{e^\prime} = \emptyset$.
\end{itemize}
The space $\UDspace{2}{G}$ is homotopy equivalent to a wedge of two circles. For example, using the Morse matching described in~\cite{farley-sabalka-2005} (with only a small change in the enumeration of vertices of the graph): If $\{3,4\}$ is the deleted edge, and the vertex $1$ is the root of the tree $T \subset G$, then $[1,2]$ is the unique critical 0-cell; the only critical 1-cells are $[\{2,4\}, 3]$ and $[\{3,4\}, 1]$; and there are no critical 2-cells. Thus, $\UDspace{2}{G} \simeq S^1 \wedge S^1$.

\begin{figure}[htb]
 \begin{center}
 \includegraphics[scale=0.75]{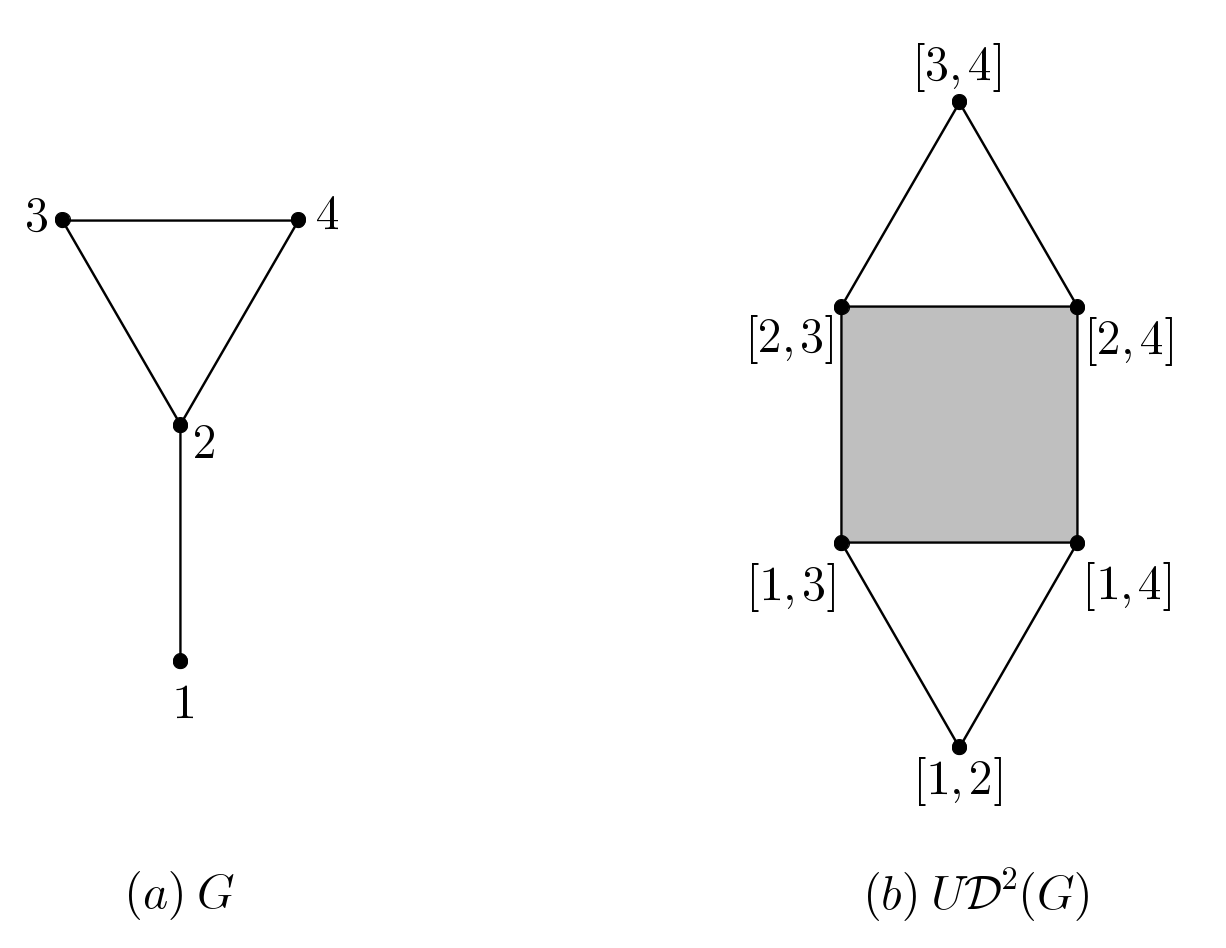}
 \end{center}
 \caption{At left is the graph $G$ from Example~\ref{collapse-example}. At right is the discretized configuration space $\UDspace{2}{G}$.}
 \label{p-optimal}
\end{figure}

Now, let $G^\prime$ be the graph obtained from $G$ by subdividing the edge $\{1,2\}$. For an explicit labeling, $G^\prime=(V,E)$, where $V=\{1, 2, 3, 4, 5\}$ and $E=\left\{ \{1,2\}, \{2,3\}, \{3,4\}, \{4,5\}, \{3,5\} \right\}$ (see Figure~\ref{p-collapse}(a)). The space $X=\UDspace{2}{G^\prime}$ is illustrated in the center, Figure~\ref{p-collapse}(b). The number of cells of each dimension in $X$ is as follows.
\begin{itemize}
 \item 0-cells: $\binom{5}{2}=10$.
 \item 1-cells: $15$.
 \item 2-cells: 4, namely, $[\{1,2\}, \{3,4\}]$, $[\{1,2\}, \{4,5\}]$, $[\{1,2\}, \{3,5\}]$, and $[\{2,3\}, \{4,5\}]$.
\end{itemize}

The gradient vector field of the discrete Morse function $f$ defined in the proof of Lemma~\ref{sclemma} is indicated by the bold arrows in Figure~\ref{p-collapse}(c). Here, the path $p \subset G^\prime$ from Lemma~\ref{sclemma} is the oriented arc from vertex $3$ to vertex $1$; vertex $2$ is the added vertex. In the edge labeling scheme from the proof of Lemma~\ref{sclemma}, $e_1 = \{2,3\}$, $\tau e_1$ is vertex $3$, $e_2 = \{1,2\}$, and $\tau e_2$ is vertex $1$.
\begin{itemize}
 \item Each arrow points from a vertex terminal cell (at its tail) to an edge terminal cell (at its head). The four arrows in Figure~\ref{p-collapse}(c) are: the upward-pointing arrow from $[1,2]$ to $[1,\{2,3\}]$, and -- from left to right -- the downward-pointing arrows from $[2,\{3,4\}]$ to $[\{1,2\}, \{3,4\}]$, from $[2,3]$ to $[3,\{1,2\}]$, and from $[2,\{3,5\}]$ to $[\{1,2\},\{3,5\}]$.
 
 \item The unmatched cells in $X$ (neither the head nor the tail of any arrow) are the critical cells of the Morse function $f$ -- the remote cells (of rank 0), which meet at most one cell of the path $p$, or meet only its boundary.
The deformation retract $Y\subset X$ consists of these remote cells, the image of the inclusion $\UDspace{2}{G} \into \UDspace{2}{G^\prime}$.
\end{itemize}

\begin{figure}[htb]
 \begin{center}
 \includegraphics[scale=0.75]{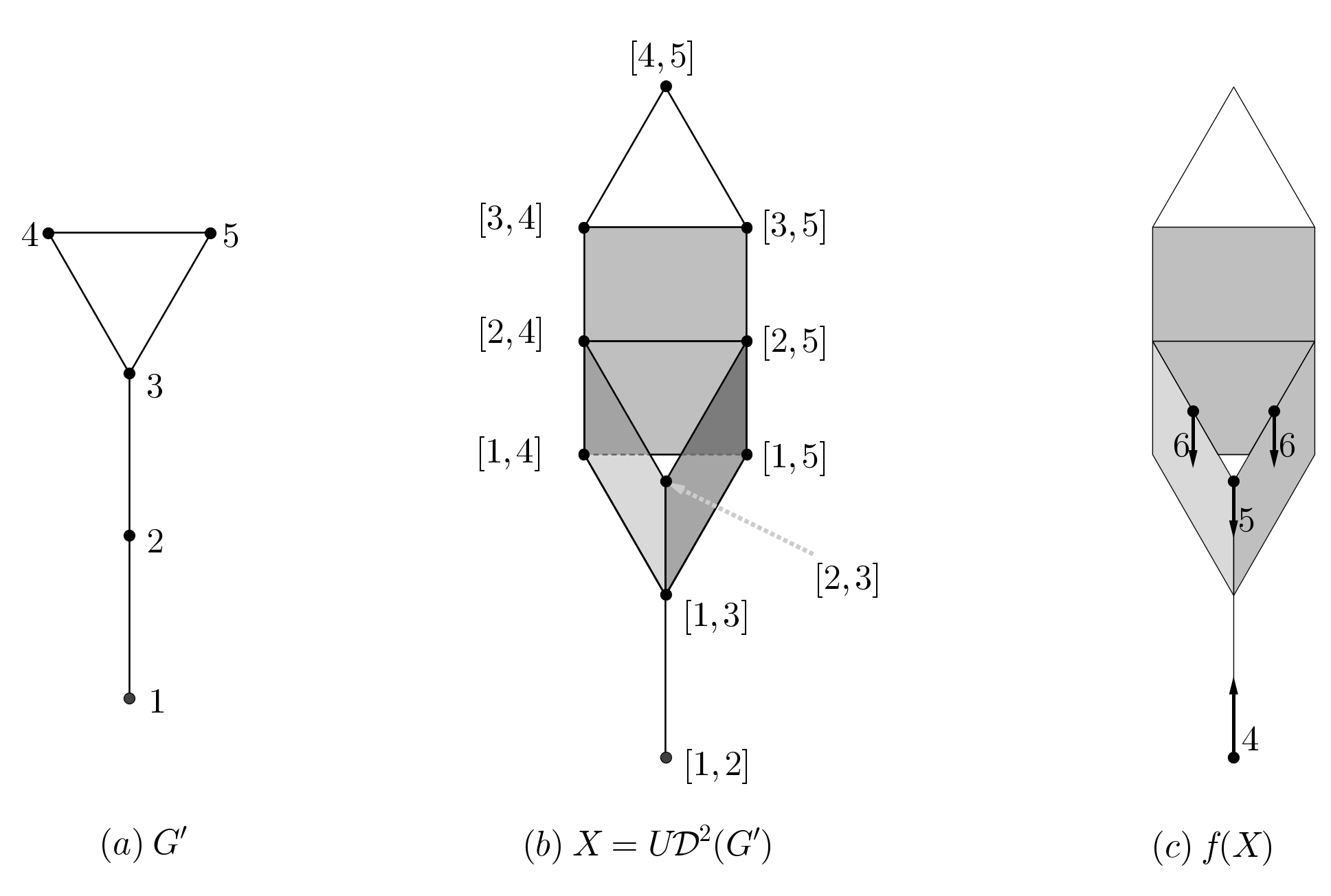}
 \end{center}
 \caption{At left (a) is the graph $G^\prime$ from Example~\ref{collapse-example}.  In the center (b) is the discretized configuration space $X = \UDspace{2}{G^\prime}$; its 0-cells are indicated by bracketed pairs $[i, j]$ of vertices. At right (c), $X$ is shown with non-critical values of $f$ and the induced Morse matching of vertex/edge terminal pairs; each critical value of $f$ is simply the dimension of the remote cell where it occurs.}
 \label{p-collapse}
\end{figure}

\label{collapse-example}
\end{example}

\bibliographystyle{plain}
\bibliography{Braidgroups}

\end{document}